\documentclass[10pt]{amsart}
\usepackage{amssymb}
\usepackage[utf8]{inputenc}
\usepackage[english]{babel}
\usepackage{latexsym,cite,mathrsfs}
\usepackage{amsfonts}
\usepackage{amssymb, amsthm}
\usepackage{xcolor}

\def\W{\accentset{\circ}{W}^{1}}
\usepackage{accents}

\numberwithin{equation}{section}

\numberwithin{equation}{section}

\usepackage{accents}

\usepackage{amssymb}
\usepackage{amsthm}

\newtheorem{thm}{Theorem}[section]
\newtheorem{lem}{Lemma}[section]

\newtheorem{prop}{Proposition}[section]

\newtheorem{rem}{Remark}[section]

\begin{document}

\title[On finding bifurcations for non-variational elliptic systems]{On finding bifurcations for non-variational elliptic systems by the extended quotients method}

\author{Yavdat ~Il'yasov}
\address{Institute of Mathematics of UFRC RAS, 112, Chernyshevsky str., 450008 Ufa, Russia}
\email{ilyasov02@gmail.com}


%

\begin{abstract}

We develop a novel method for finding bifurcations for nonlinear systems of equations based on directly finding bifurcations through saddle points of extended quotients.  The method is applied to find the saddle-node bifurcation point for elliptic equations with the nonlinearity of the general convex-concave type. The main result justifies the variational formula for the detection of the maximum saddle-node type bifurcation point of stable positive solutions.  As a consequence, a precise threshold value separating the interval of the existence of stable positive solutions is established.

\end{abstract}


\maketitle

\section{Introduction}
This paper develops a method of detecting bifurcation introduced in \cite{IlyasFunc, IlyasJDE24}, which provides a direct way of finding  bifurcations by identifying saddle points of the corresponding extended Rayleigh quotient. 
 We develop the method by finding saddle-node  bifurcation point for the following system of equations: 
\begin{equation}
	\label{p}
	\left\{ \begin{aligned}
		-\Delta u_i  &=   a_i(x) u_i^q + \lambda g_i(x,u),  &&x \in \Omega, \\[0.4em]
		~~~u_i\geq &~ 0,  &&x \in \Omega, \\[0.4em]
		~~~u_i |_{\partial \Omega} &= 0,~~i=1,\ldots, m. 
	\end{aligned}\right.
\end{equation}
Here   $q_i \in (0,1)$, $i=1,\ldots, m$, $\Omega$ is a bounded connected domain in $\mathbb{R}^d$ with  $\partial \Omega \in C^2$, $d \geq 1$, $\lambda \in \mathbb{R}$,  $u:=(u_1,\ldots, u_m)$. For $i=1,...,m$, $a_i \in L^\infty(\Omega)$, $a_i > 0$ in $\Omega$,  $a_i(\cdot)$ is a H\"older continuous  function in $\Omega$,  $g_i(\cdot,u) $ is measurable function in $\Omega$, $\forall u \in \mathbb{R}^m$, and $g_i(x,\cdot) \in C^1(\mathbb{R}^m,\mathbb{R})$. Furthermore,
	
	\par 
	$(g_1):$\, $ \exists c_0,  c_1 > 0$ and $\exists \gamma \in ( 1,+\infty)$ such that 
	$$
0\leq g_{i}(x,u) \leq 	c_0|u|+c_2|u|^{\gamma},~~ x \in \Omega,~~ u \in \mathbb{R}^m_+, ~ i=1,...,m;
	$$
	\par
	$(g_2):$\, $\exists c_2, c_3>0$ and $\exists \gamma_0 \in ( 1, \gamma)$ such that 
	$$
	 (g_{i,u_i}(x,u)u_i^2- g_i(x,u)u_i)\geq c_2  |u|^{\gamma_0+1}+c_3  |u|^{\gamma+1},~~ x \in \Omega,~~ u \in \mathbb{R}^m_+.
	$$
	Throughout this paper the summation convention is in place: we sum over any index that appears twice. A particular  example of functions $g_i$, $i=1,\ldots,m$ that meets condition $(g_1)-(g_2)$ is as follows: $g_i=b_i(x) u_i+b(x)\sum_{j=1}^m |u_j|^{\gamma-1}u_j$, $i=1,\ldots,m$ with $\gamma_0=\gamma >m$ and $b, b_i \in L^\infty(\Omega)$.

	 Hereafter, we denote  $\mathcal{W}:=(\W_2(\Omega)\cap L^\gamma)^m$, $F_i(u,\lambda):= -\Delta u_i-   a_i |u_i|^{q-1}u_i- \lambda g_i(x,u)$, $i=1,\ldots m$, $F(u,\lambda)=(F_1(u,\lambda),\ldots,F_m(u,\lambda))^T$, $u \in \mathcal{W}$, $\lambda \in \mathbb{R}$. $(\mathcal{W})^*$ means the dual space  of $\mathcal{W}$.  
A point $(u, \lambda) \in \mathcal{W}\times \mathbb{R}$ is called a weak solution  of \eqref{p} if  $F(u, \lambda)=0$ holds true in  $(\mathcal{W})^*$.

Define
$
S:=\{u \in  C^1(\overline{\Omega}) \mid \exists c_u>0,  ~u > c_u \mbox{ dist}(x,\partial \Omega)~~\mbox{in}~\Omega, ~u|_{\partial\Omega}=0\}.
$ 
We show  below that the map $F(\cdot,\lambda):\mathcal{W} \to \mathcal{W}^*$ is  Fr\'echet differentiable on $S^m$, $\forall \lambda \in \mathbb{R}$ (see Proposition \ref{prop:cont} below). A solution $u_\lambda \in S^m$ of \eqref{p} is said to be \textit{stable} if $\lambda_1(F_u(u_\lambda,\lambda))\geq 0$, and \textit{asymptotically stable} if $\lambda_1(F_u(u_\lambda,\lambda))>0$, cf. \cite{CranRibin, Dupaigne}. Hereafter,  $\lambda_1(F_u(u,\lambda))$, for $u \in S^m$ denotes the first eigenvalue of the operator $F_u(u,\lambda)(\cdot)$. 

Introduce 
\begin{align*}
	&\mathcal{W}_s:=\{u \in \mathcal{W}\cap S^m:~~\lambda_1(F_u(u,\tau))\geq  0, ~\tau=\mathcal{R}(u, u)\},
\end{align*}
We call a solution $(\hat{u},\hat{\lambda}) \in \mathcal{W}_s\times \mathbb{R}$ of \eqref{p} the \textit{saddle-node  bifurcation point}  in $\mathcal{W}_s$ (or, equivalently, fold, turning point) (cf. \cite{keller1, kielh})   if the following is fulfilled: (i) the nullspace $N(F_u(\hat{u},\hat{\lambda}))$ of the Fr\'echet derivative $F_u(\hat{u},\hat{\lambda})$ is not empty; (ii) there exists $\varepsilon >0$  and a neighborhood $U_1 \subset \mathcal{W}$  of $\hat{u}$ such that   for any $\lambda \in (\hat{\lambda},\hat{\lambda}+\varepsilon)$ equation \eqref{p} has no solutions in $\mathcal{W}_s\cap U$; (iii) for each $\lambda \in (\hat{\lambda}-\varepsilon, \hat{\lambda})$, the equation has precisely two distinct solutions in $\mathcal{W}_s\cap U$. This definition corresponds to the solution's curve turning back at the bifurcation value $\hat{\lambda}$. The solution's curve turning forward  is defined similar.  In the case only (i)-(ii) are satisfied, we call $(\hat{u},\hat{\lambda})$ the \textit{saddle-node type bifurcation point} of \eqref{p} in $\mathcal{W}_s$. 
A saddle-node type bifurcation point   $(u^*,\lambda^*)$ is said to be \textit{maximal} in $\mathcal{W}_s$ if  $\hat{\lambda}\leq  \lambda^*$ for any other saddle-node  type  bifurcation $(\hat{u},\hat{\lambda})$ of \eqref{p} in $\mathcal{W}_s$.

A model  example for \eqref{p} in the scalar case  is the Ambrosetti–Brezis–Cerami problem
\cite{ABC} with concave–convex nonlinearity  
\begin{equation}
	\label{ps}
	-\Delta u=   |u|^{q-1}u + \lambda |u|^{\gamma-1}u, ~~~u |_{\partial \Omega} = 0,
\end{equation}
where $0<q<1<\gamma$. It is why the nonlinearity in \eqref{p} can be considered to be of the convex-concave type. From \cite{ABC} it follows that there exists an extremal  value $\lambda^*>0$ such that for any $\lambda \in (0,\lambda^*]$, \eqref{ps} has a stable positive   solution $u_\lambda$, while for  $\lambda>\lambda^*$, \eqref{ps} does not admit  weak positive solutions. According to \cite{ABC}, the solution $u_\lambda$ for $\lambda \in (0,\lambda^*)$ is obtained by super-subsolution methods, while $u_{\lambda^*}$ is shown to exist as a limit point of $(u_\lambda)$. Unfortunately, this method is not easily adaptable to systems of equations like \eqref{p}. Indeed, the super-sub solution method for a system of equations differs considerably from that which is used for a scalar equation.

In general cases, system \eqref{p} is not a variational or Hamiltonian. It should be noted that in contrast to the extensive literature concerning the existence of solutions for variational and Hamiltonian systems (see the survey \cite{Figueiredo}), relatively little research is devoted to nonvariational and non-Hamiltonian systems of equations (see, e.g., \cite{alves, Caffar, Clap, crooks, ural, wu} and references therein).


The finding of bifurcations of solutions to equations poses a more complex challenge, requiring a comprehensive approach that considers both the finding solutions themselves and the analysis of the structure of the family of solutions. This problem is still quite challenging even when dealing with scalar equations. The complete answer to the question on the existence of the saddle-node bifurcation point and the exact shape of the positive solution curves for instance  of the scalar equation \eqref{ps} was obtained only in radially symmetric solutions \cite{Korman, Ouyang, Tang}.  An additional obstacle encountered when studying equations \eqref{p} and \eqref{ps} is the presence of singular  derivatives of the right-hand sides. Specifically, standard methods (see \cite{CranRibin, CranRibinExtr, keller1, kielh}) are not readily applicable for verifying that the solution's curve turning back at $(\hat{u},\hat{\lambda})$, due to the difficulties in testing conditions for the second derivative of $F(u,\lambda)$.

The existence of positive solutions and multiplicity results were studied only in some special cases of  system \eqref{p} in the variational form (see, e.g., \cite{He,Hsu,IlyasJDE24,wu} and references therein). A recent study \cite{IlyasJDE24} answered the question of whether positive solutions of  system \eqref{p} in the variational form have a saddle-node type bifurcation point. However, in the general cases of system \eqref{p}, to the best of my knowledge, no studies have been conducted on the existence of non-negative solutions and saddle-node bifurcation points. 



%
Let us state our main results. Observe that by the definition the saddle-node type bifurcation point $(u, \lambda) \in \mathcal{W}_s\times \mathbb{R}$ of \eqref{p} should  satisfy the system of equations
\begin{equation}\label{eq:Msys} 
	\left\{ \begin{aligned} 	&F(u, \lambda)   =0,\\ 	&F_u(u, \lambda) (v)=0, 	
	\end{aligned}\right.
\end{equation}
with some $v \in N(F_u(u, v) )$. To analyze this system, following \cite{IlyasJDE24} we introduce the \textit{extended Rayleigh quotient} (\textit{extended quotient} for short) associated with \eqref{p}
\begin{align*}
	\notag
	\mathcal{R}&(u, v) :=  \frac{\int (\nabla u_i, \nabla v_i) \, 
		- \int  a_i u_i^{q}v_i}
	{\int g_i(x,u)    v_i },~~u\in \mathcal{W}_s, v\in \Sigma(u).
\end{align*}
Here $\Sigma(u):=\{v \in \mathcal{W}: \int g_i(x,u)    v_i   \neq 0\}$ for $ u \in \mathcal{W}_s$.
Observe,  
\begin{align*}
	\left\{ \begin{aligned}
		&\lambda=\mathcal{R}(u,v),\\
		&\mathcal{R}_v(u,v)=0,\\
		&\mathcal{R}_u(u,v)=0
	\end{aligned}\right.~\Leftrightarrow~
	\left\{ \begin{aligned} 	&F(u, \lambda)   =0,\\ 	&F_u(u, \lambda) (v)=0, 	
	\end{aligned}\right.
\end{align*}
that is,  the set of saddle-node type bifurcation points of \eqref{p} contains in the set of critical points of $\mathcal{R}(u,v)$ on $\mathcal{W}_s\times  \mathcal{W}$.

In our approach, the following minimax formula plays a major role (cf. \cite{IlyasJDE24})
\begin{equation}\label{MainB}
	\lambda^*_s:= \sup_{u\in \mathcal{W}_s}\inf_{v\in \Sigma(u)}\mathcal{R}(u, v).
\end{equation}

The main result of the work is as follows
\begin{thm}\label{thmM} 
	Assume  $(g_1)$ - $(g_2)$, $q_i \in (0,1)$, $i=1,\ldots, m$.
	\par
	$(1^o)$  Then $0\leq  \lambda^*_s < +\infty$.  
	 	  	\begin{description}
	  		\item[{\rm (a) }] For $\lambda=\lambda^*_s$, there exists a weak positive   solution $u^*_{s} \in (C^{1,\alpha}(\overline{\Omega}))^m\cap \mathcal{W}$ of system \eqref{p}. 
			\item[{\rm (b) }] For any $\lambda>\lambda^*_s$, system \eqref{p} has no stable weak positive     solutions. 
	  	\end{description}
	
	\par 
	$(2^o)$ Assume in addition that $\gamma<2^*$, and   $q_i < (2^*-2)/2\equiv 2/(d-2)$ if $d>2$. 	
	Then $0<  \lambda^*_s < +\infty$,  $(u^*_{s}, \lambda^*_s)$ is a maximal saddle-node type bifurcation point of \eqref{p} in $\mathcal{W}_s$. Moreover, $u^*_{s}$ is a stable solution of \eqref{p}.  

		\end{thm}

Here $2^*=2d/(d-2)$ if  $d\geq 3$, and $2^*=+\infty$ if $d=1,2$.

\begin{rem} Statement $(1^o)$ can be supplemented as follows. There exists $\lambda \in [0,\lambda^*_s]$ such that \eqref{p} has a stable positive weak solution $u_\lambda \in (C^{1,\alpha}(\overline{\Omega}))^m\cap \mathcal{W}$. Indeed, we will see below that \eqref{p} has a stable positive weak solution at least for $\lambda=0$. 
\end{rem}

The following can also be considered in conjunction with the value \eqref{MainB}
\begin{equation}\label{Suppl}
	\lambda^*_{as}:= \sup_{u\in \mathcal{W}_{as}}\inf_{v\in \Sigma(u)}\mathcal{R}(u, v).
\end{equation}
Here $\mathcal{W}_{as}:=\{u \in \mathcal{W}\cap S^m:~~\lambda_1(F_u(u,\tau))>  0, ~\tau=\mathcal{R}(u, u)\}$. It easily see that $\lambda^*_{as}\leq \lambda^*_s$. For \eqref{Suppl}, it can be obtained a result similar to Theorem \ref{thmM}. In particular, we have the following
\begin{thm}\label{thm2}
	Assume  $(g_1)$ - $(g_2)$, $q_i \in (0,1)$, $i=1,\ldots, m$, and    $q_i < (2^*-2)/2\equiv 2/(d-2)$ if $d>2$. Then $0<  \lambda^*_{as} < +\infty$, and 
	
	\begin{description}
		\item[ \rm{(a)} ] For $\lambda=\lambda^*_{as}$, system \eqref{p} has  a stable   weak positive solution $u^*_{as} \in (C^{1,\alpha}(\overline{\Omega}))^m\cap \mathcal{W}$. Furthermore,  $(u^*_{as}, \lambda^*_{as})$ is a maximal saddle-node type bifurcation point of \eqref{p} in $\mathcal{W}_{as}$.
	\item[ \rm{(b)} ] For any $\lambda>\lambda^*_{as}$, \eqref{p} has no asymptotically stable weak positive solutions.  
	\item[ \rm{(c)}]  There exists a sequence of  asymptotically stable weak  positive solutions $u_{\lambda_n}\in (C^{1,\alpha}(\overline{\Omega}))^m\cap \mathcal{W}$ of \eqref{p} with $\lambda=\lambda_n>0$, $n=1,\ldots$ such that  $u_{\lambda_n} \to u^*_{as}$  in $\mathcal{W}$    and $\lambda_n \to \lambda^*_{as}$ as $n\to +\infty$. 
	\end{description}
\end{thm}

\begin{rem}
	It is natural to expect that  $u^*_s=u^*_{as}$, $\lambda^*_s= \lambda^*_{as}$ and $(u^*_s, \lambda^*_s)$ is indeed a saddle-node bifurcation point of \eqref{p} in $\mathcal{W}_s$.
	It should be noted that assertions $(1^o)$, (b) and (2) of Theorems \ref{thmM} and \ref{thm2} do not necessarily mean that \eqref{p} has no any positive solutions for $\lambda>\lambda^*_s$. Furthermore, such a behavior is possible if \eqref{p} has an $S$-shaped bifurcation curve (see \cite{BrownIbShiv,  IlCarvaSant, Gilmore}).
	
\end{rem} 
\begin{rem}
We believe that the variational formula \eqref{MainB} has the potential  provide a useful tool in further analyzing saddle-node bifurcation points and constructing numerical methods for finding them (cf. \cite{  IlD, IlIvan1, IlIvan2, ilBifChaos,  Salazar}).
\end{rem} 
The rest of the paper is organised as follows. Section \ref{sec: thm0} presents preliminaries.   
In Section 3, we prove Theorems \ref{thmM}, \ref{thm2}. In Appendix, we present a proof of a version of  Ekeland's principal for smooth functional.

\section{Preliminaries}\label{sec: thm0}
We use the standard notation $L^p:=L^p(\Omega)$
for the Lebesgue spaces, $1 \leq p \leq +\infty$, and denote by $\|\cdot\|_p$
the associated norm. By $\accentset{\circ}{W}_2^{1}:=\accentset{\circ}{W}_2^{1}(\Omega)$ we denote the standard Sobolev space, endowed with the norm $\|u \|_{1,2}=(\int |\nabla u|^2)^{1/2}$. Hereafter, we denote $W:=\W_2(\Omega)\cap L^\gamma$, $d(x):=\mbox{ dist}(x,\partial \Omega)$.  

For $\delta>0$, define  $
S(\delta):=\{u \in  S \mid~u(x) > \delta d(x)~\mbox{in}~\Omega\}$,  $S(0):=S$. Clearly,  $S(\delta)$, $\forall \delta \geq 0$ is an open subset in $C^1:=C^1(\overline{\Omega})$.
Let $Y$ be a topological space such that $S^m(\delta) \subset Y$.  The set $S^m(\delta)$ endowed with topology of $Y$ we denote by $S^m_{Y}(\delta)$. $\mathcal{L}(\mathcal{W},\mathcal{W}^*)$ denotes the Banach space of bounded linear operators from $\mathcal{W}$ into  $\mathcal{W}^*$. 

\begin{prop}\label{prop:cont} Let  $\lambda \in \mathbb{R}$, $0<q_i<1$, $i=1,\ldots,m$. 
Then $F_i(u,\lambda):\mathcal{W} \to \mathcal{W}^*$ is Fr\'echet differentiable at any $u\in S^m$, and $F_i(\cdot,\lambda) \in C^1(S^m_{C^1}, \mathcal{L}(\mathcal{W},\mathcal{W}^*))$. Furthermore, if in addition $q_i\leq (\bar{\gamma}-2)/2$, where $\bar{\gamma}=\max\{\gamma,2^*\}$, then  $F_{i,u}(\cdot,\lambda) \in C(S_{\mathcal{W}}(\delta), \mathcal{L}(\mathcal{W},\mathcal{W}^*))$, $\forall \delta>0$. 
\end{prop}
\begin{proof}  We develop an approach proposed  in \cite{ABC}. We verify the assertion only for the map $Q(u):=a_i u^q$, $u \in S$ since for the remaining terms in $F_i$ the statement is trivial. 
	Using the inequality $u(x)\geq c(u)d(x)$ in $\Omega$ and  H\"older's inequalities we derive

	\begin{align*}
		|\int a_i u^{q-1}\phi\psi  \, dx| =&|\int a_i u^q\left(\frac{\phi}{u}\right) \psi \, dx|		\leq\\
		& \frac{\|a_i\|_{\infty}}{c(u)}\|u\|^q_{p}\cdot\|\frac{\phi}{d(\cdot)}\|_2\cdot \|\psi\|_{\gamma}, ~~u \in S,~\forall \phi, \psi\in C^\infty_0(\Omega), 
	\end{align*}

	where $p=2q\gamma/(\gamma-2)$.  By the Hardy inequality, $\|\phi/d(\cdot)\|_2\leq C\|\phi\|_{1,2}$,  $\forall \phi \in C^\infty_0(\Omega)$, and thus, using the Sobolev inequalities we derive
		\begin{equation}\label{welldef}
		|\int a_i u^{q-1}\phi\psi \, dx | \leq \frac{C\|a_i\|_{\infty}}{c(u)} \|u\|^q_{p
		}\|\phi\|_{1,2}\|\psi\|_{1,2}, ~~u \in S,~\forall \phi,\psi\in C^\infty_0(\Omega),
	\end{equation}
	where $C\in (0,+\infty)$  does not depend on $u, \phi, \psi$.
	This implies that $Q(\cdot):W \to W^*$ is Fr\'echet differentiable at any $u\in S$. In the same manner we can see  that   $Q_u(\cdot) \in  C(S_{C^1}, \mathcal{L}(S_{C^1},W^*))$, and thus, $F_i(\cdot,\lambda) \in C^1(S^m_{C^1}; \mathcal{L}(\mathcal{W},\mathcal{W}^*))$.
	
Let us prove the second part. For simplicity we assume that $\bar{\gamma}=\gamma$. Suppose $u_k \to u$ in $L^\gamma$ as $k \to +\infty$. This means that  there exist $\bar{u} \in L^\gamma$ and a subsequence (which is denoted again by $(u_k)$) such that $|u_k|, |u| \leq \bar{u}$  in $\Omega$. Hence, $u_k^{q-1}d(x)\leq \bar{u}^q$  in $\Omega$, $k=1,\ldots $, and therefore,   Lebesgue's dominated convergence theorem yields $u_k^{q-1}d(x) \to u^{q-1}d(x)$ in $L^{\gamma/q}$ as $k\to +\infty$. 
	 Notice that $q\leq  ({\gamma}-2)/2$ implies $p\leq {\gamma}$. Similar to \eqref{welldef}
	we have
	\begin{equation}\label{welldef2}
	|	\int a_i (u_k^{q-1}-u^{q-1})\phi\psi \, dx| \leq \frac{C}{c(u)} \|(u_k^{q-1}-u^{q-1})d(\cdot)\|^q_p
		\|\phi\|_{1,2}\|\psi\|_{1,2}, ~~~\forall \phi,\psi\in W,
	\end{equation}
	for some $C<+\infty$ which does not depend on $u, u_k \in S$, $\phi, \psi\in (W)^*$. 	
	 Hence, using the Sobolev inequalities we derive $Q_u(\cdot)\in C(S_{W}(\delta), \mathcal{L}(\mathcal{W},\mathcal{W}^*))$, and thus,  $F_{i,u}(\cdot,\lambda) \in C(S_{\mathcal{W}}(\delta), \mathcal{L}(\mathcal{W},\mathcal{W}^*))$, $\forall \delta>0$.
\end{proof}

\begin{prop}\label{prop1}
	If 	$u \in  \mathcal{W}$ is a weak non-negative solution   to \eqref{p}, then $u_i \in C^{1,\alpha}(\overline{\Omega})$ for any $\alpha \in (0, 1)$, and $u_i>0$ in $\Omega$, $i=1,\ldots,m$. 
\end{prop} 
\begin{proof} Note that equality \eqref{p}  implies that $-\Delta u_i\geq 0$, $i=1,\ldots,m$, and thus, by the maximum principals for the elliptic problems, $u_i>0$ in $\Omega$, $i=1,\ldots,m$.
	The standard bootstrap argument and  Sobolev's embedding theorem entail that  $u_i \in L^\infty(\Omega)$, $i=1,\ldots,m$.   This by the regularity results for elliptic problems from \cite{LadUral} implies that $u_i \in C^{1,\alpha}(\overline{\Omega})$ for any $\alpha \in (0, 1)$, $i=1,\ldots,m$.
\end{proof}
Let $i=1,\ldots,m$. By Brezis-Oswald's  result \cite{BrezOsw} there exists a unique solution $w_i \in \W_2(\Omega)\cap C^{1,\alpha}(\overline{\Omega})$, $\alpha \in (0, 1)$ of 
\begin{equation}
	\label{q}
	\left\{ \begin{aligned}
		-\Delta w=& a_i w^q~~\mbox{in}~~ \Omega,\\
		w |_{\partial \Omega}&= 0.
	\end{aligned}\right.
\end{equation}
By the assumption $a_i(\cdot)$ is a H\"older continuous  function in $\Omega$, and hence, by the Schauder estimates  (see, e.g., \cite{Giltrud}), $w_i \in C^{2}(\Omega)$. Furthermore, the strong maximum principals for the elliptic problems imply that   $\displaystyle{\min_{x' \in \partial \Omega}\frac{\partial w_i(x')}{\partial \nu(x')} >0}$, where $\nu(x')$ denotes the interior unit normal at $x' \in \partial \Omega$, see, e.g., Lemma 3.4 in \cite{Giltrud}.
Thus, $w_i \in S$. 

Moreover,  $\overline{w}:=(w_1,\ldots,w_m)$ is a stable solution of \eqref{q} with $\lambda=0$, i.e., $\overline{w} \in \mathcal{W}_s$. Indeed, from Proposition \ref{prop:cont} it follows that $F(\overline{w},0) \in C^1(S^m_{C^1}; \mathcal{L}(\mathcal{W},\mathcal{W}^*))$, and therefore,  $\lambda_1(F_u(\overline{w},0))$ is well defined. It is not hard to show  that $\overline{w}$ is a minimizer of  $E(v):=\left(\frac{1}{2}\sum_{i=1}^m\int |\nabla v_i|^2-\frac{1}{q+1}\int a_i |v_i|^{q+1} \right)$ on $(\W_2)^m$, that is, 
$E(\overline{w})=\inf_{v \in (\W_2)^m}E(v)$. In particular, this means 
\begin{equation}\label{eq:lambda1}
	\lambda_1(F_u(\overline{w},0))\geq  0.
\end{equation}

\begin{lem}\label{lem1} Suppose $u^0 \in \mathcal{W}_s$ such that
	$$
	-\infty<\lambda_0:=\inf_{v \in \Sigma(u^0)}\mathcal{R}(u^0 ,v)<+\infty.
	$$
	Then  
	$u^0$ is a weak solution of \eqref{p} for $\lambda=\lambda_0$.
	
\end{lem}
\begin{proof}
	Let 	 $v^k \in \Sigma(u^0)$, $k=1,\ldots$, such that
	$$
	\lambda_k\equiv \mathcal{R}(u^0, v^k) \to  \inf_{v \in \Sigma(u^0)}\mathcal{R}(u^0,v)\equiv\lambda_0~~\mbox{as}~~ k \to +\infty.
	$$
	Since $\mathcal{R}(u, v)=\mathcal{R}(u, s v)$, $\forall s \in \mathbb{R}\setminus 0$, $\forall v \in \Sigma(u)$, $\forall u  \in \mathcal{W}_s$, we may assume that
	\begin{equation}\label{Norm}
		\int  g_i(x, u^0)v^k_i=1,~~ k=1,\ldots.
	\end{equation}
	Calculate
	\begin{align*}
		\mathcal{R}_{v}&(u^0, v^k)({\xi})=\frac{\int(\nabla u^0_i, \nabla \xi_i)-\int a_i (u_{i}^0)^{q}\xi_i-\mathcal{R}(u^0, v^k)\int g_i(x, u^0)\xi_i}{\int g_i(x, u^0)v_i^k },\\
		\mathcal{R}_{vv}&(u^0, v^k)({\xi}, {\zeta})=\\
		&-\frac{\left(\int(\nabla u^0_i, \nabla \xi_i)-\int a_i (u_{i}^0)^{q}\xi_i-\lambda_k\int g_i(x, u^0)\xi_i\right)\cdot \int g_i(x, u^0)\zeta_i}{(\int g_i(x, u^0)v_i^k)^2 }-\\
		& \frac{\left(\int(\nabla u^0_i, \nabla \zeta_i)-\int a_i (u_{i}^0)^{q}\zeta_i-\lambda_k\int g_i(x, u^0)\zeta_i\right)\cdot \int g_i(x, u^0)\xi_i}{(\int g_i(x, u^0)v_i^k )^2 },~\zeta, \xi \in \mathcal{W}.
	\end{align*}
	Let $\phi \in \mathcal{W}$, $\|\phi\|_{\mathcal{W}}=1$. Using \eqref{Norm} and the H\"older and Sobolev inequalities one can see that
	\begin{align}\label{nomina}
		|\int ( g_i(x, u^0)(v^k_i+\tau {\phi}_i)|=&|1+\tau \int g_i(x, u^0) {\phi}_i |\geq   1-a_0|\tau|,
	\end{align}
	where $ a_0 \in (0, \infty)$ does not depend on   ${\phi}$ and $k=1,\dots$. Hence   $v^k+\tau {\phi} \in \Sigma(u^0)$ for any $k=1,\dots$ and $\tau$ such that $|\tau|<\tau_0:=1/a_0$.

	By \eqref{nomina}, $(g_1)$ we have
	\begin{align}\label{DDEst}
		\|\mathcal{R}_{vv}({u}^0,  v^k+\tau \phi)\|_{(\mathcal{W}\times \mathcal{W})*}=\frac{2}{|\int g_i(x,u^0) (v^k_i+\tau {\phi}_i) |^2}\times \quad \quad \quad \quad\quad \quad& 	\\
		\sup_{\xi, \zeta \in \mathcal{W}}\frac{|\left(\int(\nabla {u}^0_i, \nabla \xi_i)-\int a_i (u_{i}^0)^{q}\xi_i-\lambda_k\int g_i(x,u^0)\xi_i\right)\cdot \int  a_i (u_{i}^0)^{q}\zeta_i|}{\|\xi_i\|_{\mathcal{W}} \|\zeta_i\|_{\mathcal{W}}}\leq &\nonumber\\
		\frac{2}{(1-a_0|\tau|)^2}\left(\sum_{i=1}^m\|-\Delta {u}^0_i -  a_i (u_{i}^0)^{q}-\lambda_kg_i(x,u^0)\|_{\mathcal{W}*} \right)\|{u}^0\|_{\mathcal{W}}&\quad \quad\quad\quad\quad\nonumber\\
		\leq \frac{C_0}{(1-a_0|\tau|)^2}, &\nonumber
	\end{align}
	where $ C_0 \in (0, \infty)$ does not depend on $k=1,\dots$. We thus may apply Theorem \ref{thm:Ek} to the  functional $G(v):=\mathcal{R}({u}^0, v)$
	defined in the open domain   $V:= \Sigma(u^0)\subset \mathcal{W}$. Indeed, it is easily seen that  $G \in C^2(\Sigma(u^0))$, and \eqref{DDEst} implies \eqref{DDRstG}, while by \eqref{nomina} there holds \eqref{DDRstG2}. Thus,  we have
	$$
	\epsilon_k:=\|\mathcal{R}_{v}(u^0, v^k)\|_{\mathcal{W}*} \to 0~~~\mbox{as}~~~k \to +\infty,
	$$
	which by \eqref{Norm} yields: 
	\begin{equation*}
		|\int(\nabla u^0_i, \nabla \xi)-\int a_i (u_{i}^0)^{q}\xi-\int \lambda_k g_i(x,u^0)\xi|\leq \epsilon_k \|\xi\|_{\mathcal{W}},~~ ~\forall \xi \in W.
	\end{equation*}
	$i=1,\ldots,m$.
	Now passing to the limit as $k \to +\infty$ we obtain \eqref{p}. 
	
\end{proof}

\section{Proof
	of main results} \label{sec: profthm1}
	
{\it Proof of Theorem \ref{thmM}}

Let us prove  $(1^o)$. Let $\overline{w}:=(w_1,\ldots,w_m)$, where $w_i$, $i=1,\ldots,m$ is a solution of \eqref{q}. By the above  $\overline{w} \in \mathcal{W}_s$, and thus, we have   
		\begin{align*}
			\lambda^*_s:= \sup_{u\in \mathcal{W}_s}\inf_{v\in \Sigma(u)}\mathcal{R}(u, v)\geq \inf_{v\in \Sigma(\overline{w})}	\frac{\int (\nabla w_i, \nabla v_i) \, 
				-  \int a_i w_i^{q}v_i}
			{\int g_i(x,\overline{w})    v_i }=0.
		\end{align*}
	Since $0\leq \lambda^*_s\leq +\infty$,  
  there exists a maximizing sequence 
		$u^n \in \mathcal{W}_s$, $n=1,\ldots$, such that 
		$$
		\lambda_n:=\lambda(u^n):=\inf_{v \in \Sigma(u^n)}\mathcal{R}(u^n,v) \to \lambda^*_s~~\mbox{as}~~n\to +\infty.
		$$ 
 By Lemma \ref{lem1}, 
	\begin{equation}\label{BEq1}
			\left\{ \begin{aligned}
			-\Delta u_i^n  &= a_i (u_i^n)^{q} + \lambda_n g_i(x,u^n)  ,  &&x \in \Omega, \\[0.4em]
			~u_i^n |_{\partial \Omega} &= 0,~~~i=1,\ldots, m. 
		\end{aligned}\right.
		\end{equation}
	 Testing  \eqref{BEq1} by $u_i$, $i=1,\ldots, m$ and integrating by parts we derive
	\begin{equation}\label{eq:oce1}
		\|u^n\|_{1,2}^2=\int a_i |u^n_i|^{q_i+1} +\lambda_n \int g_i(x,u^n) u_i^n, ~~n=1,\ldots.
	\end{equation} 
	Since $u^n \in \mathcal{W}_s$, $, ~~n=1,\ldots$,
	\begin{equation}\label{eigenProb}
		\lambda_1(F(u^n,\lambda_n)):=\inf_{\phi \in \mathcal{W}}\frac{\int |\nabla \phi|^2 -q_i \int a_i (u^n_i)^{q_i-1}|\phi_i|^2 -\lambda_n\int g_{i,u_{j}}(u^n)\phi_j\phi_i }
		{\int |\phi|^2 \, }\geq 0.
	\end{equation}
	Hence
	\begin{equation}\label{eq:oce2}
		\|u^n\|_{1,2}^2 \geq q_i \int a_i |u^n|^{q_i+1}+ \lambda_n\int g_{i,u_{i}}(u^n)(u^n_i)^2 \, ~~n=1,\ldots.
	\end{equation}
	Subtraction \eqref{eq:oce1} from \eqref{eq:oce2}, and using $(g_2)$ we obtain 

 \begin{align*}
 	(1-q_i)\int a_i |u^n|^{q_i+1}\geq \lambda_n   &\int \left( g_{i,u_{i}}(u^n)(u^n_i)^2-g_i(x,u^n) u_i^n \right)\geq\\
 	& \lambda_n (c_2 \|u^n\|_{\gamma_0}^{\gamma_0}+c_3 \|u^n\|_{\gamma}^{\gamma}), ~~n=1,\ldots. 
 \end{align*}

	Applying H\"older's inequality we derive
	\begin{equation}\label{eq:enerEst}
		C_1 \geq  \lambda_n (c_2 \|u^n\|_{\gamma_0}^{\gamma_0-q}+c_3 \|u^n\|_{\gamma}^{\gamma-q}), ~~n=1,\ldots.
	\end{equation}
	Let us show that there exists $C_2 \in (0,+\infty)$ which does not depends on $n=1,\ldots,$ such that 
	\begin{equation}\label{eq:separZERO}
		\|u^n\|_{\gamma_0}, \|u^n\|_{\gamma}\geq C_2, ~~n=1,\ldots,.
	\end{equation}

	To this end we need the following assertion that is derived by the same method as Lemma 3.3 in \cite{ABC}
	\begin{lem} \label{ABCLem}
		Assume that $f(t)$ is a function such that $t^{-1}f(t)$ is decreasing for $t>0$, $a \in L^\infty$, $a> 0$ in $\Omega$. Let $v$ and $w$ satisfy:
		$u>0$, $w>0$ in $\Omega$, $v=w=0$ on $\partial \Omega$, and
		\begin{align*}
			-\Delta w\leq a f(w),~~-\Delta u\geq a f(u), ~~\mbox{in}~~ \Omega.
		\end{align*}
		Then $u\geq w$.
	\end{lem}

	By the assumption  $g_i(x,u)\geq 0$, $x \in \Omega$, $i=1,\ldots,m$, $u \in \mathbb{R}^m_+$, and therefore,
	$$ 
	-\Delta u_i^n  \geq    a_i (u_i^n)^q ~~  in ~~ \Omega,~~i=1,\ldots,m,~~ n=1\ldots.
	$$
	Hence,  \eqref{q} and Lemma  \ref{ABCLem} yield
	\begin{equation}\label{ineqABC}
		u_i^n \geq w_i, \quad i=1,\ldots,m,~~ n=1\ldots, 
	\end{equation}
and as a result, we get \eqref{eq:separZERO}. Clearly, \eqref{eq:enerEst}, \eqref{eq:separZERO} imply that $\lambda^*_s<+\infty$, and  $ \|u^n\|_{\gamma}\leq C_2$, $n=1,\ldots$. This by \eqref{eq:oce1}, $(g_1)$ implies that  $\|u^n\|_{1,2}\leq C_3$, $n=1,\ldots$, where $ C_2, C_3 \in (0,+\infty)$ do not depend on $n=1,\ldots$. Thus, $(u^n)$ is bounded in $ \mathcal{W}$, and  therefore, by the Banach–Alaoglu and the Sobolev theorems  there exists a subsequence (again denoted by $(u^n)$) such that
	\begin{align}
		&	u^n \rightharpoondown u^*_{s}~~\mbox{weakly in $\mathcal{W}$},   \label{wse}
		\\		&u^n \to u^*_{s}~~\mbox{strongly in}~~ (L^r)^m, ~~1\leq r< 2^*, \label{sse} 
	\end{align} 
$\mbox{as }~~n\to +\infty$ for some $ u^*_{s}\in \mathcal{W}$. 	From \eqref{ineqABC} it follows tha $u^*_{s,i}\geq w_i>0$, $i=1,\ldots,m,~~ n=1\ldots$.
	
Passing to the limit in \eqref{BEq1} as $n\to +\infty$ we obtain
	\begin{equation}\label{eq:fequat}
		-\Delta u^*_{s,i}  =  a_i (u_i^*)^{q} + \lambda^*_sg_i(x, u^*_{s}),~~x \in \Omega,~~i=1,\ldots,m.
	\end{equation}
Using Proposition \ref{prop1} we conclude that   $u^*_{s} \in C^{1,\alpha}(\overline{\Omega})$ for any $\alpha \in (0, 1)$ and $u_i^*>0$ in $\Omega$, $i=1,\ldots,m$. Thus,  we  obtain $(1^o)$, (a).

To show  $(1^o)$, (b), suppose conversely that for $\lambda>\lambda^*_s$ there exists a stable weak non-negative solution $u_\lambda$ of \eqref{p}. Then by Proposition \ref{prop1}, $u_{\lambda}\in \mathcal{W}_s$, and consequently,  \eqref{MainB} yields $\inf_{v\in \Sigma(u_\lambda)}\mathcal{R}(u_\lambda, v)<\lambda$. Hence, there exists $v \in \Sigma(u_\lambda)$ such that $\mathcal{R}(u_\lambda, v)<\lambda$. Assume that  $\int u_{\lambda,i}^{q}v_i>0$. Then
$$
\int(\nabla u_{\lambda,i}, \nabla v_i)-\int a_i u_{\lambda,i}^{q}v_i-\lambda\int g_i(x, u_\lambda) v_i<0,
$$
which contradicts \eqref{p}, and as a result we get   $(1^o)$, (b).

Let us prove $(2^o)$. For simplicity we assume that $d>2$.   Using  $\gamma < 2^*$, it is not hard to show from \eqref{wse},\eqref{sse}, \eqref{eq:fequat}, and $u^*_{s}\neq 0$  that 
\begin{equation}\label{eq:strConv}
	u^n \to u^*_{s}~~\mbox{ strongly in}~~ \mathcal{W}~\mbox{as}~n\to +\infty.
\end{equation}
Clearly, by the  maximum principle  $w_i \in S(\delta)$, $i=1,\ldots,m$ with some $\delta>0$. Hence,  \eqref{ineqABC} imply that $u_n \in S^m(\delta)$, $n=1,\ldots $. 
The assumption    $q_i < (2^*-2)/2\equiv 2/(d-2)$ for $d>2$ and  $\gamma < 2^*$ implies by Proposition \ref{prop:cont} that $F_u(u,\lambda) \in C(S^m_{\mathcal{W}}(\delta); \mathcal{L}(\mathcal{W},\mathcal{W}^*))$. Hence, 
\begin{align}\label{ConvDR}
	\langle F_u(u_n,\lambda_n)(\phi),\psi \rangle \to
	\langle F_u(u^*_{s},\lambda^*_s)(\phi),\psi \rangle~ \mbox{as}~~ n\to +\infty&,\quad\forall \phi,\psi  \in \mathcal{W},
\end{align}
and consequently, $\lambda_1(F_u(u_n,\lambda_n)) \to \lambda_1(F_u(u^*_{s},\lambda^*_s))\geq 0$ as $n\to +\infty$. Thus we get that $u^*_{s} \in \mathcal{W}_s$.

Let us show that 
$\lambda_1(F_u(u^*_{s},\lambda^*_s))=0$.
Suppose, contrary to our claim, that $\lambda_1(F_u(u^*_{s},\lambda^*_s))>0$. Then    $F_u(u^*_{s},\lambda^*_s)(\cdot): \mathcal{W} \to \mathcal{W}^*$  is nonsingular linear operator. From Proposition \ref{prop:cont} we have  $F(\cdot,\lambda) \in C^1(S^m_{C^1};\mathcal{L}(\mathcal{W},\mathcal{W}^*))$. 
Hence, by the Implicit Functional Theorem (see, e.g., \cite{Dieudonn}) there is a neighbourhood $V\times U \subset \mathbb{R}\times S^m_{C^1}$  of $(\lambda^*_s, u^*_{s})$ and a mapping $V \ni \lambda \mapsto u_\lambda \in U$ 
such that $u_\lambda|_{\lambda=\lambda^*_s}=u^*_{s}$ and $F(u_{\lambda}, \lambda)=0$, $\forall \lambda \in V$. 
Furthermore, the map $u_{(\cdot)}: V \to U$ is continuous. Since $\lambda_1(F_u(u^*_{s},\lambda^*_s))>0$, there exists a neighbourhood $V_1 \subset V$ of $\lambda^*_s$ such that $\lambda_1(F_u(u_{\lambda},\lambda))>0$ for every $\lambda \in V_1$. However, this contradicts assertion ($1^o)$ of the theorem. Thus,  $\lambda_1(F_u(u^*_{s},\lambda^*_s))=0$, and $(u^*_{s},\lambda^*_s)$ is a saddle-node type bifurcation point of \eqref{p} in $\mathcal{W}_s$. Since $(1^o)$, (b), $(u^*_{s},\lambda^*_s)$ is a maximal saddle-node type bifurcation point of \eqref{p} in $\mathcal{W}_s$.

Finally, let us show that 	$0<\lambda^*_s$. Suppose the converse $\lambda^*_s=0$. Then by the above $0=\lambda_1(F_u(u^*_{s},\lambda^*_s))=\lambda_1(-\Delta-q|u^*_{s}|^{q-1})$. However, this is clearly impossible.

\medskip 

{\it Proof of Theorem \ref{thm2}} is similar to the proof  of Theorem \ref{thmM}. We only need to show assertion (c) of Theorem. Indeed, by  the construction there are sequences $\lambda_n$, $u^n_{as}  \in \mathcal{W}_{as}$, $n=1,\ldots$ such that  $F_u(u^n_{as},\lambda_n)=0$, $n=1,\ldots$, $\lambda_n \to \lambda^*_{as}$,  
and $u^n_{as} \to u^*_{as}$  strongly in $\mathcal{W}$ as $n\to +\infty$. Moreover, $u^*_{as} \in \overline{\mathcal{W}}_{as}$, and $\lambda_1(F_u(u^*_{as},\lambda^*_{as}))=0$. Thus, $u^*_{as}  \in \overline{\mathcal{W}}_{as} \setminus \mathcal{W}_{as}$. On the other hand,  
$u^n_{as}  \in \mathcal{W}_{as}$, $n=1,\ldots$. Hence, $u^n_{as}  \neq u^*_{as} $, $n=1,\ldots$, and we thus obtain the proof of (c).

\section{Appendix A}\label{sec:appendix}
Let $X$ be a Banach space and  $V\subset X$ be an open  domain. Denote  $B_r:=\{\phi \in X: ~~\|\phi\|_{X}\leq r\}$, $r>0$. Assume that $G: V \to \mathbb{R}$, $G \in C^2(V)$. 
Consider
\begin{equation}\label{ApMin}
	\hat{G}=\inf_{v \in V}G(v).
\end{equation}

\begin{thm}\label{thm:Ek}
	Assume that $|	\hat{G}|<+\infty$. Suppose that there exist $\tau_0, a_0, C_0\in (0,+\infty)$, and a minimizing sequence $(v_k) \subset V$ of \eqref{ApMin} such that  
	\begin{align}
		&\|G_{vv}(v_k+\tau \phi)\|_{(X\times X)^*} < 	\frac{C_0}{(1-|\tau| a_0)^2}<+\infty,\label{DDRstG}\\
		&v_k+\tau \phi \in V,~~ \forall \tau \in (-\tau_0,\tau_0), ~~\forall \phi \in B_1,~~\forall k=1,\ldots.\label{DDRstG2}
	\end{align}
	Then 
	$$
	\|G_{v}(v_k)\|_{X^*}	:=\sup_{\xi \in X\setminus 0}\frac{|G_{v}(v_k)(\xi)|}{\|\xi\|_X} \to 0~~~\mbox{as}~~~ k \to +\infty.
	$$
\end{thm}
\begin{proof} 
	Suppose, contrary to our claim, that there exists $\alpha>0$ such that
	$$
	\|G_{v}(v_k)\|_{X^*}>\alpha,~~~\forall k=1,\ldots. 
	$$
	This means that for every $k=1,\ldots$, there exists $\phi_k \in V$, $\|\phi_k\|_{X}=1$ such that
	$|G_{v}(v_k)(\phi_k)|>\alpha$. 
	By the Taylor expansion  
	$$
	G(v_k +\tau \phi_k)=G(v_k)+\tau G_{v}(v_k)(\phi_k)+\frac{\tau^2}{2}G_{vv}(v_k+\theta_k \tau \phi_k)(\phi_k,\phi_k),
	$$
	for sufficiently small $|\tau|$, and some $\theta_k \in (0,1)$, $k=1,\ldots$. Suppose, for definiteness, that  $G_{v}(v_k)(\phi_k)>\alpha$. Then for $\tau\in (-\tau_0, 0)$, by \eqref{DDRstG}
	\begin{align*}
		G(v_k +\tau \phi_k)\leq G(v_k)+\tau \alpha +\frac{\tau^2}{2}\frac{C_0}{(1+\tau a_0)^2}, ~~~k=1,\ldots. 
	\end{align*}
	It is easily seen that there exists $\tau_{1}\in (0,\tau_0)$  such that  
	$$
	\kappa(\tau):=\tau \left( \alpha +\frac{\tau}{2}\frac{C_0}{(1+\tau a_0)^2}\right) <0,~~~~\forall \tau \in (-\tau_{1}, 0).
	$$
	Since $(v_k)$ is a minimizing sequence, for any $\varepsilon>0$ there exists $k(\varepsilon)$ such that 
	$$
	G(v_k)<\hat{G}+\varepsilon,~~\forall k>k(\varepsilon).
	$$
	Take $\tau \in (-\tau_1, 0)$ and  $\varepsilon_0=-\kappa(\tau)/2$. Then by the above
	$$
	G(v_k +\tau \phi_k)<\hat{G}+\varepsilon_0+\kappa(\tau) = \hat{G}+\kappa(\tau)/2<\hat{G},~~\forall k>k(\varepsilon_0),
	$$
	and thus, in view of \eqref{DDRstG2} we get a contradiction.
\end{proof}

\end{document}